\documentclass[12pt]{article}

\usepackage{hyperref}
\usepackage{graphicx}
\usepackage{amssymb}
\usepackage{amsfonts}
\usepackage{moreverb}

\def\complib{COMP{l$_\mathrm{e}$ib}}
\def\matlab{{\sc matlab}}
\def\hifoo{{\sc hifoo}}

\begin{document}

\begin{center}
\Large\bf Some control design experiments with \hifoo\
\end{center}

\begin{center}
Didier Henrion$^{1,2}$
\end{center}

\footnotetext[1]{LAAS-CNRS, 7~Avenue du Colonel Roche, 31077~Toulouse,
France. E-mail: {\tt henrion@laas.fr}}

\footnotetext[2]{Also with the Faculty of Electrical Engineering,
Czech Technical University in Prague, Czech Republic}

\begin{center}
\today
\end{center}

\section{Introduction}

In \cite{hifoo} a new \matlab\ package called \hifoo\
was proposed for $H_{\infty}$ fixed-order controller design.
This document illustrates how some standard controller design examples
can be solved with this software.

\section{Static output feedback design}

Consider plant {\tt AC1} from the \complib\ database \cite{complib}, a
fifth-order system modeling an aircraft transfer function with three
inputs and three outputs. Let us first stabilize this plant by
static output feedback:
\begin{verbatim}
% retrieve (A,B,C) data from COMPLIB
>> [A,dummy,B,dummy,C] = compleib('AC1');
% build (A,B,C,D) plant with D=0
>> P = ss(A,B,C,zeros(size(C,1),size(B,2)))
a = 
            x1       x2       x3       x4       x5
   x1        0        0    1.132        0       -1
   x2        0  -0.0538  -0.1712        0   0.0705
   x3        0        0        0        1        0
   x4        0   0.0485        0  -0.8556   -1.013
   x5        0  -0.2909        0    1.053  -0.6859
b = 
            u1       u2       u3
   x1        0        0        0
   x2    -0.12        1        0
   x3        0        0        0
   x4    4.419        0   -1.665
   x5    1.575        0  -0.0732
c = 
       x1  x2  x3  x4  x5
   y1   1   0   0   0   0
   y2   0   1   0   0   0
   y3   0   0   1   0   0
d = 
       u1  u2  u3
   y1   0   0   0
   y2   0   0   0
   y3   0   0   0
Continuous-time model.
% compute stabilizing SOF
>> K = hifoo(P,'+')
hifoo: found a stabilizing controller , quitting
d = 
             u1        u2        u3
   y1    0.1778  -0.06802     -2.76
   y2    0.6741    -1.402     2.051
   y3     1.463     2.957    -1.568
Static gain.
% check closed-loop stability
>> T = feedback(P,-K); % positive feedback
>> eig(T)
ans =
  -0.2537 + 3.2758i
  -0.2537 - 3.2758i
  -2.3229          
  -0.0796 + 1.1206i
  -0.0796 - 1.1206i
% this is the same as eig(A+B*K.d*C)
\end{verbatim}
Note that the stabilizing SOF gain computed on your own
computer generally differs from the one above, since the optimization
algorithms in \hifoo\ are randomized.

\section{Spectral abscissa minimization for the two-mass-spring
problem}

We consider a system
consisting of two masses interconnected
by a spring, a typical control benchmark problem which is a generic
model of a system with a rigid body mode and one vibration mode
\cite{aiaa}.
If the first mass is pulled sufficiently far apart from the second mass and
suddenly dropped, then the two masses will oscillate until they reach
their equilibrium position.

The control problem consists of appropriately moving the second mass
so that the first mass settles down to its final position as fast as
possible; more specifically, we want to maximize
the asymptotic decay rate. For this we use a linear feedback
controller between the system output (measured position of the second
mass) and the system input (actuator positioning the first mass).
The open-loop transfer function between the system input
and output is given by
\[
P(s) = \frac{1}{s^4+2s^2}.
\]
when mass weights and the spring constant are normalized to one.

In \cite{massspring} we show that the asymptotic decay rate
maximization problem has a non-trivial solution if and only if the
controller order is equal to two. The problem can be solved with
\hifoo\ as follows:
\begin{verbatim}
% open-loop transfer function in state-space form
>> P = ss(tf(1,[1 0 2 0 0]));
% call HIFOO
>> K = hifoo(P,2,'s');
>> tf(K)
\end{verbatim}
We find the controller
\[
\frac{6.8308175s^2-1.8486865s-0.28043397}
{s^2+4.2752492s+6.0786141}.
\]
The closed-loop poles, obtained with the
commands
\begin{verbatim}
>> T = feedback(P,-K);
>> eig(T)
\end{verbatim}
are placed at $-0.7073\pm i0.2979$, $-0.7073\pm
i0.2980$, $-0.7231\pm i0.5343$ so the achieved spectral abscissa
is $-0.7073$. We observe the typical eigenvalue clustering
phenomenon characteristic of the neighborhood of a local minimizer of
the spectral abscissa. Note that we display the
controller coefficients to 8 significant digits because
clustered eigenvalues are typically very sensitive to perturbations.
In other words, controllers obtained by optimizing the spectral 
abscissa are typically quite non-robust.

We can call \hifoo\ again with the controller just obtained as an initial
guess:
\begin{verbatim}
>> K = hifoo(P,2,'s',K);
\end{verbatim}
This yields another controller
with an improved closed-loop spectral abscissa of $-0.7380$.
One more run of \hifoo\ produces a controller
\[
\frac{8.073790s^2-1.7330367s-0.23544720}
{s^2+4.5435259s+6.7343390}
\]
further pushing the spectral abscissa to $-0.7572$.

As explained in \cite{massspring},
solving the pole placement equation by hand, assigning all the poles to
the same negative real number $-\alpha$ (a unique pole of multiplicity
six), we obtain analytically $\alpha=\frac{\sqrt{15}}{5} \approx
0.7746$ and the locally optimal second-order controller
\[
\frac{\frac{43}{5}s^2-\frac{54\sqrt{15}}{125}s-\frac{27}{125}}
{s^2+\frac{6\sqrt{15}}{5}s+7}
\approx
\frac{8.6000s^2-1.6731s-0.2160}{s^2+4.6476s+7}.
\]
We can see that \hifoo\ found numerically a very similar controller
and that the achieved spectral abscissa was not far from the one
derived analytically.

\section{NASA HIMAT aircraft}

This example is aimed at illustrating that \hifoo\ can be used
to design robust MIMO controllers of low order ensuring
a similar performance than full order controllers designed
with standard tools.

In the user's guide of the Robust Control Toolbox for \matlab\
\cite{rct} we can find
an $H_{\infty}$ design example for controlling the pitch axis of a
NASA HiMAT aircraft.
For this plant, a 10th-order controller can be designed with the
Robust Control Toolbox as follows:
\begin{verbatim}
% open loop plant
>> ag =[-2.2567e-02 -3.6617e+01 -1.8897e+01 -3.2090e+01 3.2509e+00 -7.6257e-01;
9.2572e-05 -1.8997e+00 9.8312e-01 -7.2562e-04 -1.7080e-01 -4.9652e-03;
1.2338e-02 1.1720e+01 -2.6316e+00 8.7582e-04 -3.1604e+01 2.2396e+01;
0 0 1.0000e+00 0 0 0;
0 0 0 0 -3.0000e+01 0;
0 0 0 0 0 -3.0000e+01];
>> bg = [0 0; 0 0; 0 0; 0 0; 30 0; 0 30];
>> cg = [0 1 0 0 0 0; 0 0 0 1 0 0];
>> dg = [0 0; 0 0];
>> G = ss(ag,bg,cg,dg);
% weight for sensitivity function = low-pass 0.15 rad/s
>> MS = 2; AS = .03; WS = 5;
>> s = tf([1 0],1); W1 = (s/MS+WS)/(s+AS*WS);
% weight for complementary sensitivity function = high-pass 400 rad/s
>> MT = 2; AT = .05; WT = 20;
>> W3 = (s+WT/MT)/(AT*s+WT)
% full-order (10th order) controller designed with Robust Control Toolbox
>> [K10,CL10,perf10] = mixsyn(G,W1,[],W3);
>> perf10
    0.7885
% Validation in time-domain and frequency-domain
>> L = G*K10; I = eye(size(L));
>> S = feedback(I,L); T = I-S;
>> figure; step(T,2);
>> figure; subplot(1,2,1); sigma(S); subplot(1,2,2); sigma(T)
\end{verbatim}
\begin{figure}
\begin{center}
\includegraphics{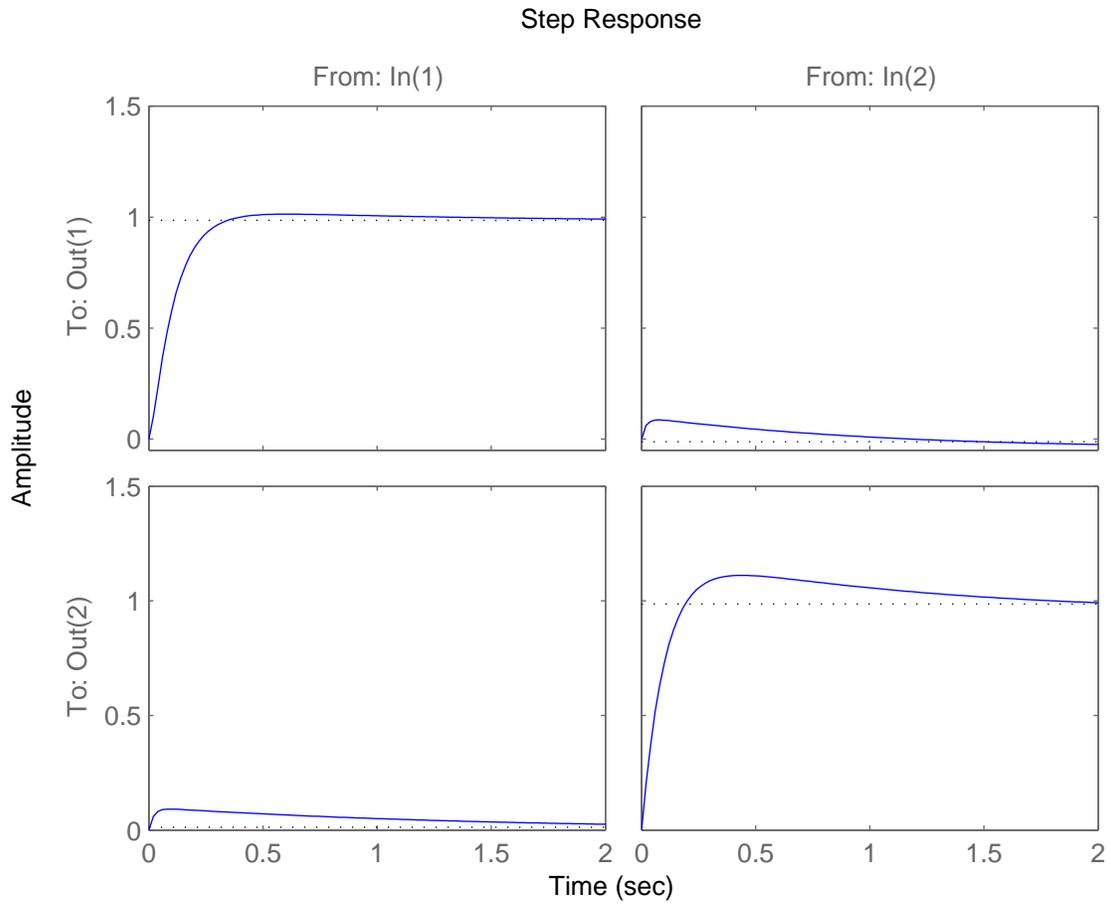}
\caption{NASA HiMAT aircraft model: step responses with the 10th order
controller designed with the Robust Control Toolbox.\label{nasa-fig1-10-fig}}
\end{center}
\end{figure}
\begin{figure}
\begin{center}
\includegraphics{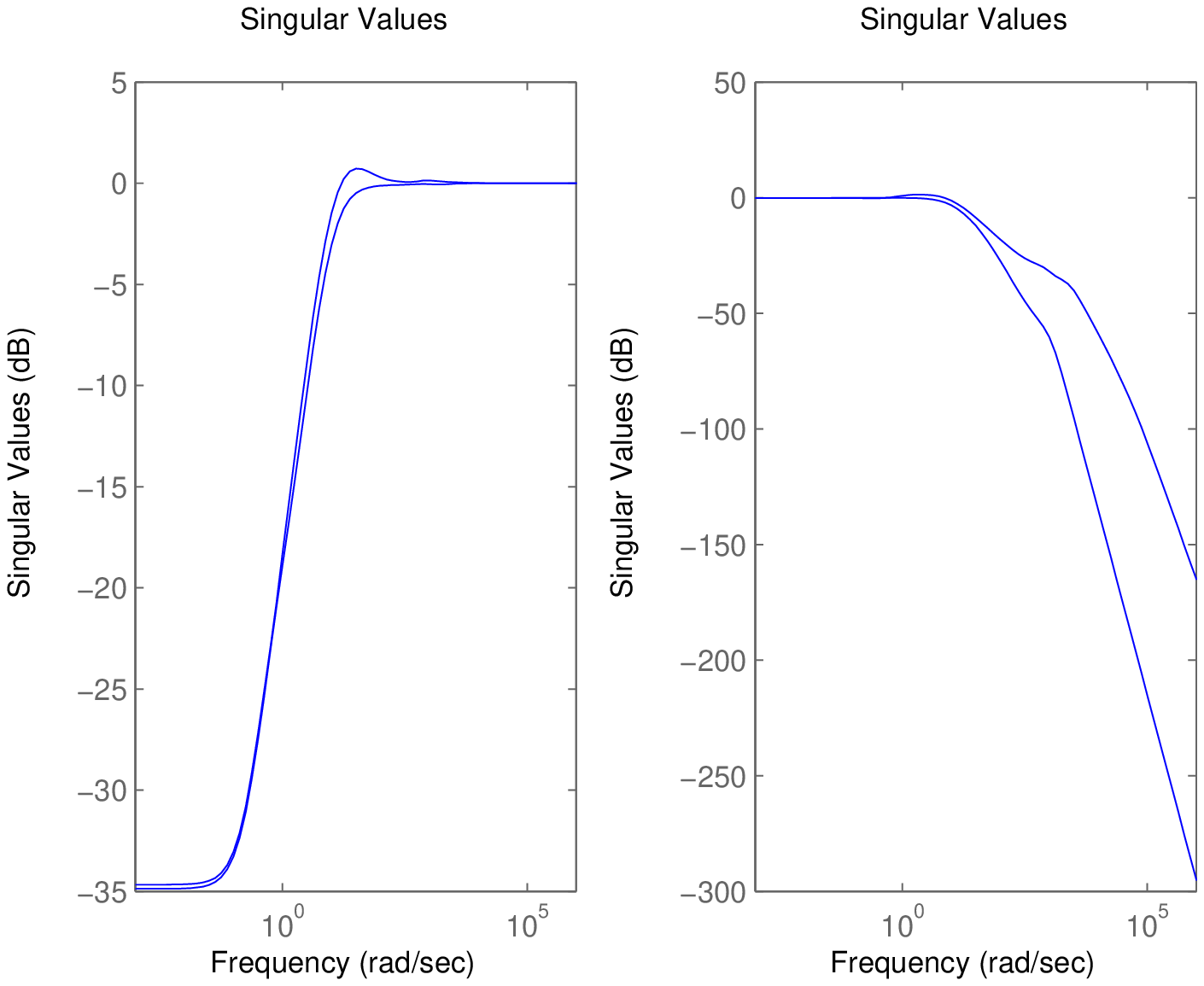}
\caption{NASA HiMAT aircraft model: singular value plots of the
  sensitivity function (left) and complementary sensitivity function
  (right) with the 10th order controller designed with the Robust
  Control Toolbox.\label{nasa-fig2-10-fig}}
\end{center}
\end{figure}
We can observe on Figure \ref{nasa-fig1-10-fig} that the step
responses are well decoupled, fast enough (raise time around 1.5
second) and without oscillations or too much overshoot/undershoot.
On Figure \ref{nasa-fig2-10-fig}, we see the typically high-pass
profile of sensitivity function $S$ (left), and the low-pass profile
of complementary sensitivity function $T$ (right), with almost no resonance.

With \hifoo\, the design can be carried out as follows:
\begin{verbatim}
% build augmented plant
>> P = augw(G,W1,[],W3);
% call HIFOO to design zero-order controller
>> [K0,perf0] = hifoo(P);
>> perf0
    3.8207
% not really satisfying, so let us try a first-order controller
>> [K1,perf1] = hifoo(P,1);
>> [K1,perf1] = hifoo(P,1,K1);
% ... after some iterations ..
>> perf1
    2.1477
% not really satisfying either..
% second-order controller
>> [K2,perf2] = hifoo(P,2);
>> [K2,perf2] = hifoo(P,2,K2);
% ... after some iterations ..
>> perf2
perf2 =
    1.5245
>> loc2
loc2 =
  2.7498e-004
% third-order controller
>> [K3,perf3] = hifoo(P,3);
>> [K3,perf3] = hifoo(P,3,K3);
% ... after some iterations ..
>> perf3
    0.9897
% Here is the third order controller:
>> K3
a = 
            x1       x2       x3
   x1    -11.1  -0.2587    31.93
   x2      2.4  0.03315   -7.116
   x3      189    2.964   -559.5
b = 
           u1      u2
   x1  -9.617   50.87
   x2   2.369  -10.98
   x3   108.1  -853.5
c = 
           x1      x2      x3
   y1   56.08   1.175  -88.97
   y2   22.51   2.271   47.12
d = 
           u1      u2
   y1  -51.53  -77.27
   y2  -106.1   156.2
Continuous-time model.
% Validation in time-domain and frequency-domain
>> L = G*K3; I = eye(size(L));
>> S = feedback(I,L); T = I-S;
>> figure; step(T,2);
>> figure; subplot(1,2,1); sigma(S); subplot(1,2,2); sigma(T)
\end{verbatim}
\begin{figure}
\begin{center}
\includegraphics{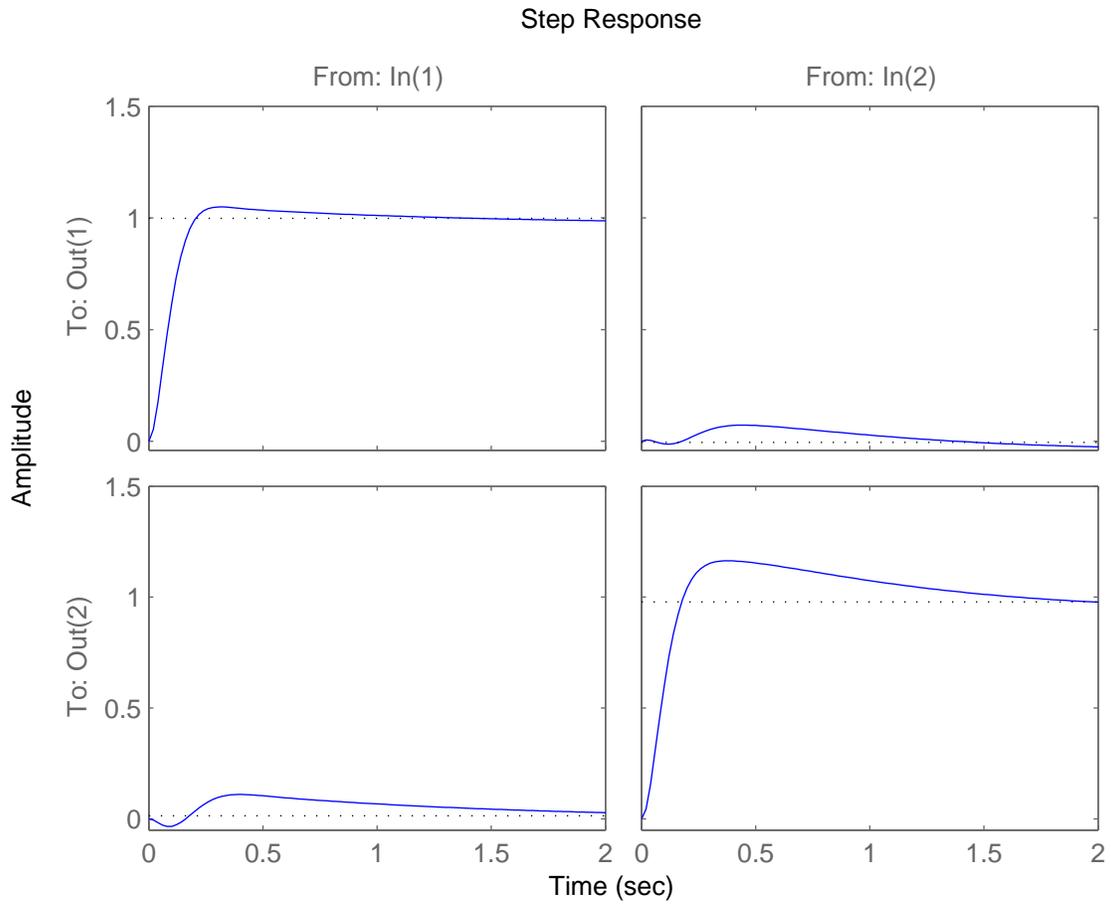}
\caption{NASA HiMAT aircraft model: step responses with the 3rd 
order controller designed with \hifoo.\label{nasa-fig1-3-fig}}
\end{center}
\end{figure}
\begin{figure}
\begin{center}
\includegraphics{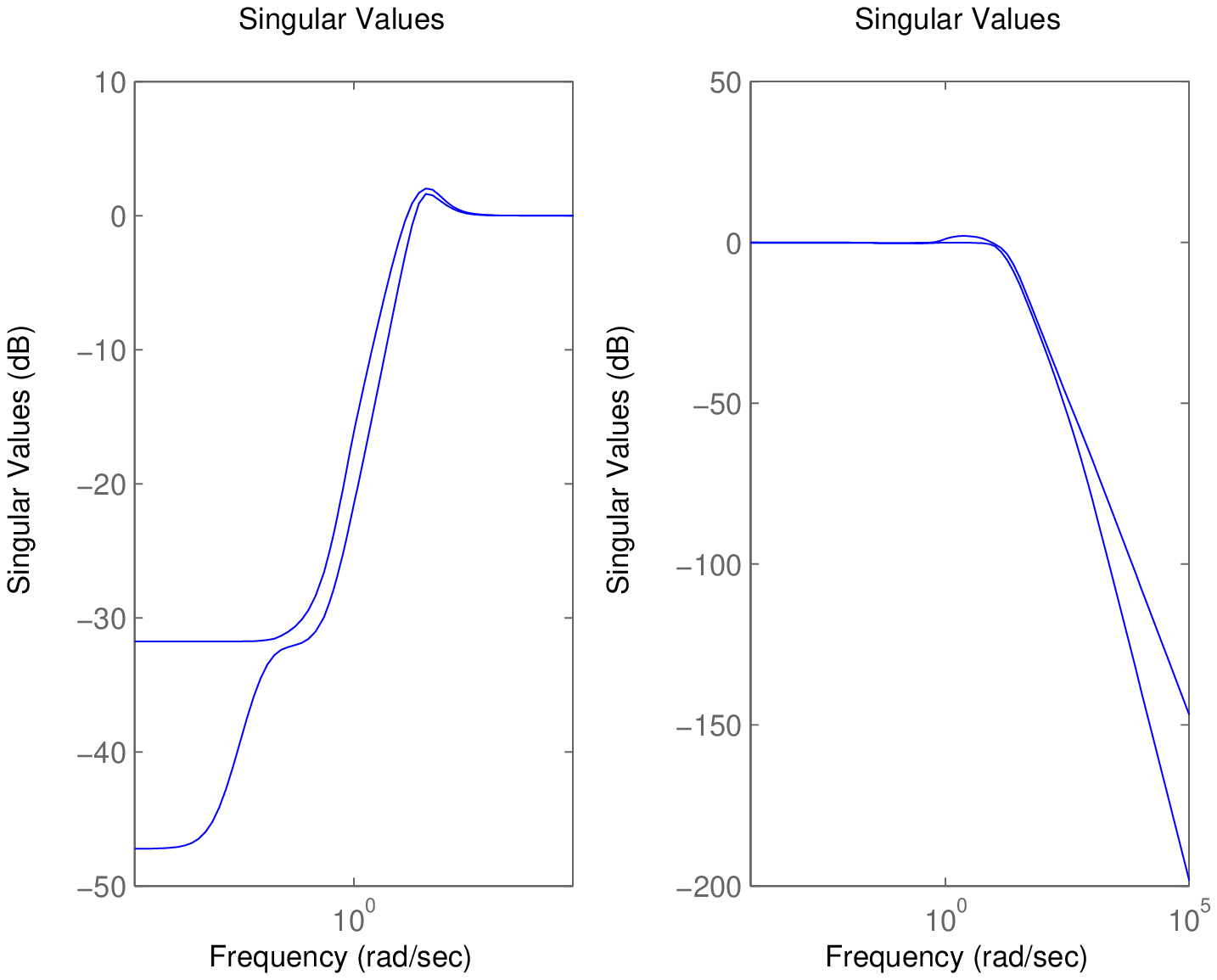}
\caption{NASA HiMAT aircraft model: singular value plots of the
  sensitivity function (left) and complementary sensitivity function
  (right) with the 3rd order controller designed with
  \hifoo.\label{nasa-fig2-3-fig}}
\end{center}
\end{figure}
Even though the $H_{\infty}$ norm achieved with the 3rd order
controller (0.9897) is still far from the norm achieved with the
10th order controller (0.7885), we can observe on Figures
\ref{nasa-fig1-3-fig} and \ref{nasa-fig2-3-fig} that the closed-loop
performance of the 3rd order controller is almost as good as that
of the 10th order controller.

\section{Four disks}

This example shows that \hifoo\ can design controllers
ensuring $H_{\infty}$ norm performances much better than
existing controller order reduction techniques.
It also illustrates the well-known (stable) pole/zero
cancellation phenomenon typical when controlling
systems with flexible modes (i.e. highly oscillatory
systems).

We consider the four-disk control system described in
\cite[Example 19.4 p. 517]{zdg}:
\begin{verbatim}
>> A = [-0.161 -6.004 -0.58215 -9.9835 -0.40727 -3.982 0 0;
eye(7) zeros(7,1)];
>> B = [1 0 1; zeros(7,3)];
>> C = [1e-3*[0 0 0 0 0.55 11 1.32 18]; zeros(1,8);
0 0 6.4432e-3 2.3196e-3 7.1252e-2 1.0002 0.10455 0.99551];
>> D = [0 0 0;0 0 1;0 1 0];
>> P = mktito(ss(A,B,C,D),1,1);
\end{verbatim}

For this plant, the best closed-loop $H_{\infty}$ norm reported
in \cite{zdg} is $1.1272$, achieved for an 8th order controller.
On pages 518-519-520 in \cite{zdg}, various controller order reduction
techniques are then applied, yielding the following $H_{\infty}$
closed-loop performances:
\begin{center}
\begin{tabular}{ll}
order & $H_{\infty}$ norm \\ \hline
8 & 1.1272 \\
7 & 1.1960 \\
6 & 1.1950 \\
2 & 1.4150 \\
1 & 2.4670 
\end{tabular}
\end{center}
With \hifoo\ we could design reduced-order controllers
ensuring the following $H_{\infty}$ performances:
\begin{center}
\begin{tabular}{ll}
order & $H_{\infty}$ norm \\ \hline
8 & 1.1317 \\
7 & 1.1267 \\
6 & 1.1326 \\
2 & 1.2438 \\
1 & 1.4256 
\end{tabular}
\end{center}
Note the significant improvement achieved for
low-order controllers. For illustration here are the first order
and second order controllers:
\begin{verbatim}
>> K1 = hifoo(P,1);
hifoo: best order 1 controller found has H-infinity performance 1.42558
>> tf(K1)
Transfer function:
-0.1227 s - 0.003706
--------------------
     s + 0.2082
>> K2 = hifoo(P,2);
hifoo: best order 2 controller found has H-infinity performance 1.24382
>> tf(K2)
Transfer function:
-0.03473 s^2 - 0.1821 s - 0.006087
----------------------------------
     s^2 + 0.6846 s + 0.2454
\end{verbatim}

\begin{figure}
\begin{center}
\includegraphics{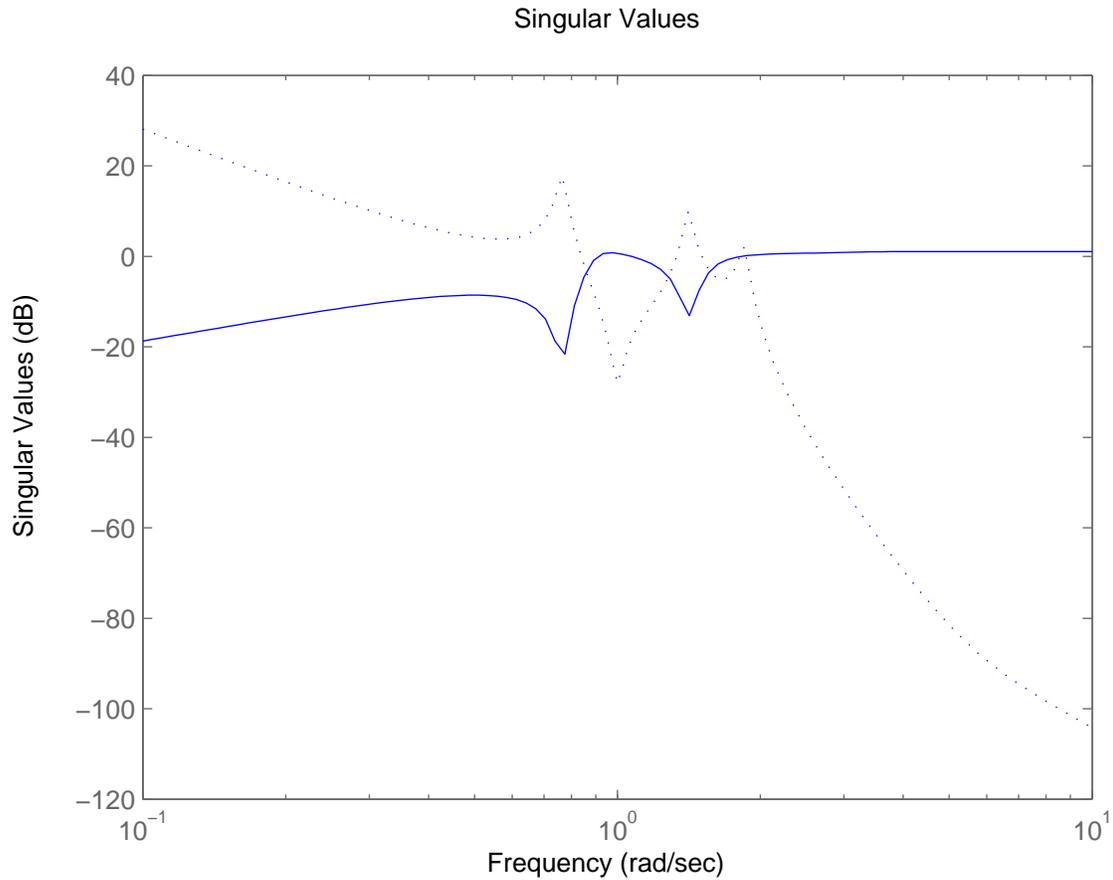}
\caption{Four-disk system: singular value plots of the 8th order
  controller (solid) and open-loop plant transfer functions
  (dotted).\label{four-disk-sigma-fig}}
\end{center}
\end{figure}

When designing a controller of sufficiently high order
we can observe the typical phenomenon of cancellation
by the controller of open-loop plant flexible modes:
\begin{verbatim}
>> K8 = hifoo(P,8);
hifoo: best order 8 controller found has H-infinity performance 1.13171
>> zpk(K8)
Zero/pole/gain:
-1.1301 (s+3.232) (s+0.03049) (s^2  + 0.02897s + 0.5845)
(s^2  + 0.08555s + 1.995) (s^2  + 2.208s + 10.51)  
--------------------------------------------------------
(s+3.295) (s+0.6869) (s^2  + 0.2009s + 0.7842)
(s^2  + 0.4285s + 2.09) (s^2  + 2.205s + 10.71)
>> zpk(P(3,3))
Zero/pole/gain:
0.0064432 (s+4.84) (s^2  + 0.04s + 1) (s^2  - 4.52s + 31.92)
------------------------------------------------------------
s^2 (s^2  + 0.0306s + 0.5852) (s^2  + 0.0564s + 1.988)
(s^2  + 0.074s + 3.423)
>> figure; sigma(K8,logspace(-1,1,100)); hold on; sigma(P(3,3),':');
\end{verbatim}
We can see that second-degree factors of the controller numerator
almost cancel second-degree factors of the plant denominator.
Singular value plots of Figure \ref{four-disk-sigma-fig}
show that indeed ``valleys'' of the controller transfer function
correspond to ``peaks'' of the open-loop plant transfer function,
meaning that the flexible modes of the plant are well damped by the
controller.
\begin{verbatim}
\end{verbatim}

\section{Order drop for optimal $H_{\infty}$ controllers}

This example illustrates a well-known feature of
$H_{\infty}$ controllers at the optimum: for regular problems,
the order of the controller is strictly less than
the order of the open-loop plant \cite{kwak93,zdg}.

When asking \hifoo\ to find a full-order controller near
the optimum, the coefficients of the controller blow up
because there is an almost pole/zero cancellation at
infinity. For ARE or LMI solvers this can be troublesome
as studied by Pascal Gahinet in the 1990s \cite{gahinet}.
Apparently \hifoo\ is fairly robust in this situation.

Consider \cite[Example 6.1]{gahinet}, a standard regular
$H_{\infty}$ problem. The open-loop plant has 3rd order,
and the $H_{\infty}$-optimal controller degenerates
to a 2nd order controller achieving an $H_{\infty}$
performance of $21.5279$. Let us try to illustrate
this phenomenon with \hifoo:
\begin{verbatim}
>> A = [1 -1 0;1 1 -1;0 1 -2]; B1 = [1 2 0;0 -1 0;1 1 0]; B2 = [1;0;1];
>> C1 = [0 0 0;1 1 0;-1 0 1]; D11 = zeros(3); D12 = [1;0;0];
>> C2 = [0 -1 1]; D21 = [0 0 1]; D22 = 0;
>> P = mktito(ss(A,[B1 B2],[C1;C2],[D11 D12;D21 D22]),size(C2,1),size(B2,2));
% Find a third-order controller
>> K3 = hifoo(P,3);
hifoo: best order 3 controller found has H-infinity performance 21.9398
>>  tf(K3)
Transfer function:
12.67 s^3 + 504.1 s^2 + 430.7 s - 632.7
---------------------------------------
   s^3 + 42.13 s^2 + 680.1 s + 64.37
>> zpk(K3)
Zero/pole/gain:
12.6721 (s+38.88) (s+1.674) (s-0.7671)
--------------------------------------
 (s+0.09521) (s^2  + 42.03s + 676.1)
>> K3=hifoo(P,3,K3);
hifoo: best order 3 controller found has H-infinity performance 21.5488
>> tf(K3)
Transfer function:
20.64 s^3 + 1097 s^2 + 961.6 s - 1362
-------------------------------------
   s^3 + 78.46 s^2 + 1475 s + 141
>> zpk(K3)
Zero/pole/gain:
20.6447 (s+52.22) (s+1.672) (s-0.7557)
--------------------------------------
   (s+47.42) (s+30.94) (s+0.09609)
\end{verbatim}
We can observe that coefficients of the transfer function $K_3$ grow in
magnitude and that there is almost a pole/zero cancellation around $s=-50$.
The theory shows that when the achieved $H_{\infty}$ norm tends to the
optimum, this pole/zero pair tends to $-\infty$ and the controller
order drops.

Let us try to design a controller with the {\tt hinfsyn} function
of the Robust Control Toolbox for \matlab, which by default solves
coupled algebraic Riccati equations (AREs) to design
full-order $H_{\infty}$-optimal controllers:
\begin{verbatim}
>> [KF,Q,perfF] = hinfsyn(P);
>> perfF % achieved H-inf norm
perfF =
   21.5284
>> tf(KF)
Transfer function:
   1.097e006 s^2 + 1.006e006 s - 1.385e006
---------------------------------------------
s^3 + 5.097e004 s^2 + 1.497e006 s + 1.435e005
>> zpk(KF)
Zero/pole/gain:
1097211.3816 (s+1.672) (s-0.7552)
----------------------------------
(s+5.094e004) (s+29.3) (s+0.09614)
\end{verbatim}
Observe how the finite zeros $(-1.7, 0.76)$ and
finite poles $(-30, -0.096)$ of $K_F$ agree with those
of $K_3$ obtained by \hifoo. Function {\tt hinfsyn}
returned a strictly proper controller $K_F$
whose frequency characteristics agree with those of $K_3$.

Now let us see what can be achieved with a second-order controller
\begin{verbatim}
>> K2 = hifoo(P,2);
hifoo: best order 2 controller found has H-infinity performance 21.5448
% Almost as good as the third order controller..
>> K2 = hifoo(P,2,K2);
hifoo: best order 2 controller found has H-infinity performance 21.5284
% ...almost reaching the optimal H-inf norm 21.5279
>> zpk(K2)
Zero/pole/gain:
21.5284 (s+1.672) (s-0.7551)
----------------------------
   (s+29.28) (s+0.09616)
\end{verbatim}
Compare with $K_3$, the poles $(-1.7, 0.76)$ and
zeros $(-30, -0.096)$ nicely fit. Note that the LMI method implemented
in function {\tt hinfsyn}, as an alternative to the ARE method,
also returns a second order controller, but the achieved performance
is a little bit worse:
\begin{verbatim}
>> [KL,Q,perfL] = hinfsyn(P,'method','lmi');
>> zpk(KL)
Zero/pole/gain:
21.5209 (s+1.672) (s-0.7543)
----------------------------
   (s+29.27) (s+0.09658)
>> perfL % Achieved H-inf norm
perfL =
  21.6040
\end{verbatim}

\section{Non-proper $H_{\infty}$ optimal controller}

This example illustrates a typical phenomenon occuring when optimizing
the $H_{\infty}$ norm of the closed-loop sensitivity function only:
the optimal controller is non-proper. This can cause numerical
troubles to ARE or LMI $H_{\infty}$ solvers, but \hifoo\ seems to be
quite robust in this case.

Consider the minimum sensitivity problem studied in \cite{kwak90},
for which we know that the optimal $H_{\infty}$ norm of $6$ is
achieved by a non-proper controller
\[
K(s) = 5 - \frac{5}{6}s.
\]
\begin{verbatim}
>> G = tf([1 -1],conv([1 -2],[1 -3]))
Transfer function:
    s - 1
-------------
s^2 - 5 s + 6
>> P = augw(G,1,[],[]);
>> K = hifoo(P,1);
hifoo: best order 1 controller found has H-infinity performance 6.01608
>> tf(K)
Transfer function:
-4.882e004 s + 2.939e005
------------------------
     s + 5.874e004
\end{verbatim}
We can see that the achieved $H_{\infty}$ is not far from the global
optimum, yet controller coefficients are very large. To understand
why, let us have a look at poles and zeros:
\begin{verbatim}
>> zpk(K)
Zero/pole/gain:
-48816.5465 (s-6.019)
---------------------
    (s+5.874e004)
\end{verbatim}
We see that there is a very large (negative) pole. We also observe
that the zero is not far from $6$. After a few more calls to HIFOO we
get something even closer to 6:
\begin{verbatim}
>> K=hifoo(P,1,K);
hifoo: best order 1 controller found has H-infinity performance 6.00024
>> zpk(K)
Zero/pole/gain:
-48748.0539 (s-6.001)
---------------------
    (s+5.851e004)
\end{verbatim}
There is still a very large pole, and a zero closer to 6.
Note that the above controller can also be written
\[
\frac{0.8332 (6.001-s)}{(1.7091\cdot 10^{-5}s+1)} \approx 
0.8333 (6.000-s)
\]
which is almost the expected globally optimal controller.

Since \hifoo\ restricts the search to controllers of first
order with monic denominator, the coefficients grow arbitrarily large
because of a pole tending to $-\infty$. It could be possible to
allow for non-proper controllers, but this is currently not
implemented in \hifoo.

\section*{Acknowledgments}

This work benefited from many interactions
with James Burke, Adrian Lewis, Marc Millstone and Michael Overton.
It was partly funded by project No.~102/06/0652 of the Grant Agency of
the Czech Republic.

\end{document}